\documentclass{article}
\usepackage{tikz}
\usepackage{amsmath}
\usepackage{amsfonts}
\usepackage{amssymb}
\usepackage{xspace}
\usepackage{graphicx}
\numberwithin{equation}{section}
\makeatletter
\newenvironment{prooof}{\par\noindent{\sc Proof:}
}{\hfill\llap{$\Box$}\vspace{1\baselineskip}\par\noindent}

\newenvironment{proofof}{\par\noindent{\sc Proof}
}{\hfill\llap{$\Box$}\vspace{1\baselineskip}\par\noindent}

\newcommand{\adl}{\vspace{1\baselineskip}}

\newtheorem{theorem}{Theorem}[section]
\newtheorem{proposition}[theorem]{Proposition}
\newtheorem{lemma}[theorem]{Lemma}
\newtheorem{corollary}[theorem]{Corollary}
\newtheorem{remark}[theorem]{Remark}
\newtheorem{definition}[theorem]{Definition}
\newtheorem{example}[theorem]{Example}

\newcommand{\beq}{\begin{equation}}
\newcommand{\eeq}{\end{equation}}

\newcommand{\ba}{\begin{array}}
\newcommand{\ea}{\end{array}}
\newcommand{\bt}{\begin{theorem}}
\newcommand{\et}{\end{theorem}}
\newcommand{\bp}{\begin{proposition}}
\newcommand{\ep}{\end{proposition}}
\newcommand{\bl}{\begin{lemma}}
\newcommand{\el}{\end{lemma}}
\newcommand{\bc}{\begin{corollary}}
\newcommand{\ec}{\end{corollary}}
\newcommand{\bi}{\begin{itemize}}
\newcommand{\ei}{\end{itemize}}
\newcommand{\ben}{\begin{enumerate}}
\newcommand{\een}{\end{enumerate}}
\newcommand{\bpf}{\begin{prooof}}
\newcommand{\epf}{\end{prooof}}
\newcommand{\bpff}{\begin{proofof}}
\newcommand{\epff}{\end{proofof}}
\newcommand{\bdf}{\begin{definition}\rm}
\newcommand{\edf}{\end{definition}}
\newcommand{\br}{\begin{remark}\rm}
\newcommand{\er}{\end{remark}}
\newcommand{\bex}{\begin{example}\rm}
\newcommand{\eex}{\end{example}}
\def\pri{\hbox to 10pt{\hfil\hbox to 0.4pt{\vrule height5pt width0.4pt
                 depth0pt}\vrule width5pt height0.4pt depth0pt\hfil}}
\newcommand{\TC}{{\rm TC}}

\newcommand{\gk}{{\bf{k}}}

\newcommand{\gK}{{\bf K}}

\newcommand{\calL}{{\cal L}}

\newcommand{\gR}{{\mathbb R}}

\newcommand{\Nat}{{\mathbb N}}

\newcommand{\RN}{{\mathbb R}^{N+1}}

\newcommand{\D}{{\cal D}}
\newcommand{\EE}{{\cal E}}
\newcommand{\EEE}{{\textrm E}}
\newcommand{\F}{{\cal F}}
\newcommand{\G}{{\cal G}}

\newcommand{\BV}{\mathop{\rm BV}\nolimits}

\newcommand{\Su}{{{\mathbb S}^1}}
\newcommand{\SN}{{{\mathbb S}^{N}}}

\def\mesh{\mathop{\rm mesh}\nolimits}

\newcommand{\kk}{{\mathfrak K}}

\newcommand{\gt}{{\bf t}}

\newcommand{\gv}{{\bf v}}

\newcommand{\gc}{{\bf c}}

\newcommand{\Var}{\mathop{\rm Var}\nolimits}

\let\a=\alpha
\let\be=\beta

\let\d=\delta
\let\e=\varepsilon

\let\vf=\varphi
\let\g=\gamma

\let\l=\lambda
\let\m=\mu

\let\p=\pi

\let\s=\sigma
\let\t=\theta
\let\tt=\tau

\let\y=\eta

\let\vf=\varphi

\let\GG=\Gamma

\let\wid=\widetilde

\let\pa=\partial
\let\sb=\subset

\let\fa=\forall

\let\sm=\setminus
\let\ol=\overline

\let\ds=\displaystyle

\let\i=\infty

\title{\Large \bf Weak elastic energy of irregular curves}
\author{\it Domenico Mucci and Alberto Saracco
\footnote{%
{\sc Dipartimento di Scienze Matematiche,
Fisiche ed Informatiche, Universit\`{a} di Parma,
Parco Area delle Scienze 53/A, I-43124 Parma, Italy.
E-mail: domenico.mucci@unipr.it, alberto.saracco@unipr.it}
}
}
\begin{document}
\date{23/12/2022}

\topskip=1.5truecm \maketitle \topskip=1.5truecm \maketitle
\maketitle
       %%%%%%%%% Abstract %%%%%%%%%%%%%%%%
%
{\small{\bf Abstract.} A weak notion of elastic energy for (not necessarily regular) rectifiable curves in any space dimension is proposed.
Our $p$-energy is defined through a relaxation process, where a suitable $p$-rotation of inscribed polygonals is adopted. The discrete $p$-rotation we choose has a geometric flavor: a polygonal is viewed as an approximation to a smooth curve and hence its discrete curvature is spread out into a smooth density.
For any exponent $p$ greater than one, the $p$-energy is finite if and only if the arc-length parameterization of the curve has a second order summability with the same growth exponent. In that case, moreover, the energy agrees with the natural extension of the integral of the $p$-th power of the scalar curvature. Finally, a comparison with other definitions of discrete curvatures is discussed.
}
\adl\par\noindent
\adl\par\noindent

%%%%%%%%%%%%%%%%%%%%

\section{Introduction}
The role played by analysis and geometry in continuum mechanics is well highlighted in the foundational work by Leonhard Euler of 1744 concerning the most classical variational model of inextensible flexible
rods.
In fact, in his Additamentum I to the monograph ``ad Methodus inveniendi lineas curvas maximi
minimive proprietate gaudentes", Euler \cite{E} stated:
 ``\textit{For a pertinent understanding of the older, fundamental works on elastics,
it is necessary to know the connections of the statements contained in them with
the methods of the Mechanics of Solids and the Mechanics of Continua.}''

\par Euler's elastica problem is briefly discussed in Sec.~\ref{Subs:Euler}, and we refer e.g. to the treatises \cite{GH,Lo,Tr} for more details.
The literature on the subject is huge and we are aware that it is not possible to give a satisfactory complete reference.
To this purpose, we address to \cite{DKS} for a study of the evolution problem of elastic curves in $\RN$, and to \cite{Miu} for
an analysis of the straightening problem, in terms of a perturbation theory for the modified total squared curvature energy.

Finally, the relaxation problem for the energy $\int_{\pa E} (1+\kk^p)\,ds$ among bounded planar sets $E$ is treated e.g. in \cite{BDMP,BM}.

\smallskip\par The aim of this paper is to give a contribution toward this direction, by proposing a weak notion of {\em bending energy}
for a wide class of {\em irregular curves}. Namely, by means of a relaxation method, for any real exponent $p>1$ we introduce a {\em $p$-curvature} functional on the class of rectifiable curves $\gc$ in $\RN$,
%denoted by $\F_p(\gc)$,
that turns out to be finite in presence of the expected {\em Sobolev regularity}, and that in the smooth case agrees with the integral of the $p$-th power of the {\em scalar curvature} $\kk$ of the curve, i.e., with the {\em $p$-energy} (or bending energy, for $p=2$)
$$ \EE_p(\gc):=\int_\gc \kk^p\,ds\,,\qquad p>1\,. $$
\par Referring to Secs.~\ref{Subs:tv}, \ref{Subs:l}, and \ref{Subs:TC} for the notation about {\em total variation}, {\em length}, and {\em total curvature}, respectively,
we remark that in the ``plastic'' case $p=1$, our functional agrees with the total curvature $\TC(\gc)$ introduced by Milnor \cite{Mi},
that is the supremum of the {\em rotation} $\gk^*(P)$, i.e., the sum of external angles (or, better, {\em turning angles})
computed among the polygonals $P$ of $\RN$ {\em inscribed} in $\gc$, say $P\ll \gc$.
\par More precisely, if $\gc$ is a rectifiable and open curve in $\RN$, we let $\gc:\ol I_L\to\RN$ denote its arc-length parameterization, so that $I_L=(0,L)$ and $L=\calL(\gc)$, the length of $\gc$. Since $\gc$ is a Lipschitz function, by Rademacher's theorem (cf. \cite[Thm.~2.14]{AFP}) it is differentiable $\calL^1$-a.e. in $I_L$, where $\calL^1$ is Lebesgue measure, so that the tangent indicatrix (or {\em tantrix}) $\gt=\dot\gc$ exists $\calL^1$-a.e. in $I_L$. If in addition $\gc$ has finite total curvature, then its tantrix is a function of bounded variation, taking values in the Gauss hyper-sphere $\SN$, and the total curvature agrees with the essential variation of $\gt$ when computed w.r.t. the intrinsic distance in $\SN$, see \eqref{VarSPt}, namely:
$$ \TC(\gc)=\Var_{\SN}(\gt)\,. $$
\par In the same spirit as Lebesgue-Serrin's relaxed functional, we introduce the {\em $p$-curvature} functional
$$ \F_p(\gc):=\inf\Bigl\{\liminf_{h\to\i}\gk_p(P_h)\mid \{P_h\}\ll\gc\,,\,\,\m_\gc(P_h)\to 0\Bigr\} \quad p\geq 1$$
of rectifiable curves $\gc$ in $\RN$, in any co-dimension $N\geq 1$, see Sec.~\ref{Subs:pcurv}.
\par In the latter centered formula, we make use of the notion by Alexandrov-Reshetnyak \cite{AR} of {\em modulus} $\m_\gc(P)$ of a polygonal $P$ inscribed in $\gc$, that is equal to the maximum of the diameter of the arcs of $\gc$ determined by the consecutive vertices in $P$.
Therefore, when dealing with the ``elastic'' case $p>1$, the first problem comes from the choice of a good notion of {\em $p$-rotation} $\gk_p(P)$
of a polygonal.
%
%\par In view of applications, in fact, the $p$-curvature of a curve can be estimated in terms of inscribed polygonals.

In discrete geometry, several definitions are proposed in order to give a discrete analogous to the $p$-energy functional $\EE_p(\gc)$ of smooth curves, see Sec.~\ref{Subs:dg}.
Taking for simplicity a closed equilateral polygonal $P$ with $n$ segments equal in length to $\ell$, and denoting by $\theta_i$ the {turning angle} at the $i$-th vertex $v_i$ of $P$,
one may e.g. choose
$$ %\beq\label{Kstarintro}
\gk_p^*(P):=\sum_{i=1}^n \frac{\theta_i^p}{(\ell/2)^{p-1}}\,,\quad p\geq 1 $$ %\eeq
so that for $p=1$ one recovers the notion of rotation $\gk^*(P)$, and hence we clearly get
$$ \F_1(\gc)=\TC(\gc)\,. $$
\par
\begin{figure}[!h]
\centering\includegraphics[width=2.5in]{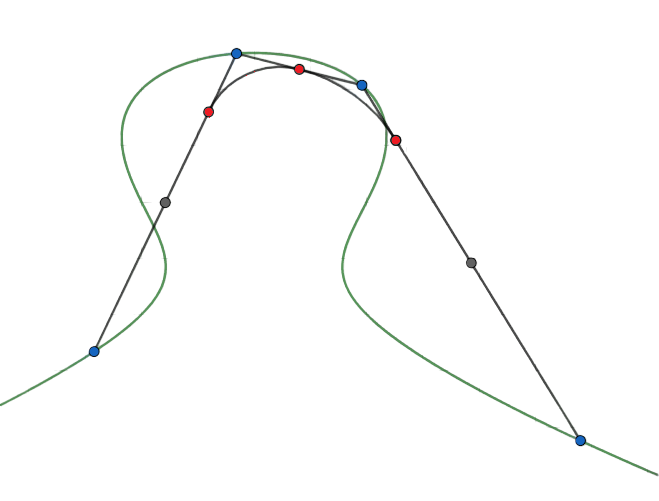}
%%% where xxxxxx name represents "figurename.eps"
\caption{A curve (green), some (blue) points on it and the poligonal given by the points. In grey the middle points of the curve and in red where the clamping condition is applied (where the distance from the blue point to the adiacent grey points is minimal). The arcs of circle smoothing the poligonal are shown.}
\label{fig:pol}
\end{figure}

 Following J. M. Sullivan \cite[Sec.~9]{Su_curv}, one may wish to view a polygonal as an approximation to a smooth curve and hence spread this curvature out into a smooth density (see Figure \ref{fig:pol}).
However, this is not the case of the discrete curvature $\gk_p^*(P)$. More precisely, when $p>1$ it is not possible to find a piecewise smooth curve $\gc$ satisfying
(first order) {\em clamping conditions} at the middle points of the edges of $P$ in such a way that its total curvature is equal to the total curvature of $P$, and its $p$-energy is equal to $\gk_p^*(P)$, see Sec.~\ref{Subs:rmks}.

\par As a simple smoothing, we may replace a neighborhood of each vertex of an equilateral polygonal (with side length $\ell$) with an inscribed circular arc touching the middle points of the two consecutive segments concurring at $v_i$. This arc turns a total angle $\t_i$, but its curvature density is $(2/\ell)\cdot\tan(\theta_i/2)$.
If e.g. $P_{n,\ell}$ is a regular $n$-agon in the plane with edges of length $\ell$, and turning angles $\t_i=2\pi/n$, for $i=1,\ldots,n$,
then its inscribed circle has curvature density equal to $(2/\ell)\cdot\tan(\pi/n)$.
\smallskip\par In this paper, we propose to define the {\em $p$-rotation} $\gk_p(P)$ of a polygonal $P$ as the {$p$-energy} of a suitable piecewise
smooth curve $\g(P)$ inscribed in $P$ by means of a generalization of the previous approach.
More precisely, the scalar curvature of $\g(P)$ turns out to be piecewise constant, since $\g(P)$ is a piecewise smooth curve given by the union of
circular arcs connected by segments, see Definition~\ref{DgP}.
Again, the curve $\g(P)$ satisfies a clamping condition at the middle points of the segments of $P$, and its total curvature agrees with the total curvature of $P$,
see Sec.~\ref{Subs:prot}. If $P$ has a turning angle equal to $\p$, we let $\gk_p(P):=+\infty$ for each $p>1$.
\par To clarify the geometric features of our construction of $p$-rotation of polygonals, some strictly related minimum problems of the $p$-energy functional $\EE_p(\gc)$ are discussed in Secs.~\ref{Subs:exmin} and \ref{Subs:exmin2}.
\par By the very definition, moreover, one readily obtains that if a rectifiable $\gc$ in $\RN$ satisfies $\F_p(\gc)<\i$ for some $p>1$, then $\F_q(\gc)<\i$ for all exponents $1\leq q<p$, and
$$ \F_q(\gc)\leq\calL(\gc)+\F_p(\gc)\,. $$
In particular, $\gc$ has finite total curvature, $\TC(\gc)<\i$, see Proposition~\ref{Ppq}.
\smallskip\par From an analytical point of view, differently from the case $p=1$, when $p>1$ a second order Sobolev regularity is expected when $\F_p(\gc)<\i$.
For this purpose, in Sec.~\ref{Subs:relcart} we collect the main results concerning a corresponding relaxed energy functional for Cartesian curves in $\RN$ that was analyzed in \cite{AcMu}, for the ``plastic case" $p=1$,
and in \cite{AcMu2}, for the ``elastic case" $p>1$.
In fact, in codimension $N=1$, as already observed by Dal~Maso et al. in \cite{DFLM}, when $p>1$, a Cartesian curve with finite relaxed ``elastic'' energy cannot have corner points.
\par In a similar way, we shall prove that a rectifiable curve has finite $p$-curvature for some $p>1$ if and only if its tantrix $\gt$ is a Sobolev map in
$W^{1,p}(I_L,\SN)$. More precisely, with the previous notation, the Main Result of this paper is enclosed in the following
\bt\label{Tmain} Let $N\geq 1$ integer, and let $\gc$ be a rectifiable and open curve in $\RN$ parameterized in arc-length. Then for every exponent $p>1$
$$ \F_p(\gc)<\i \iff \gc\in W^{2,p}(I_L,\RN) $$
see Definition~$\ref{DEp}$, and in that case
$$ \F_p(\gc)=\int_0^{L}\Vert\ddot \gc(s)\Vert^p\,ds\,. $$ \et
\par Our Main Result is coherent with the physical interpretation: an elastic rod
needs infinite bending energy in order to produce {a corner.}
Moreover, in case of smooth curves we get
$$ \F_p(\gc)=\EE_p(\gc):=\int_\gc \kk^p\,ds\,,\qquad \fa\,p>1\,. $$
\par The energy lower bound and Sobolev regularity will be proved in Sec.~\ref{Subs:lbd}, and the energy upper bound in Sec.~\ref{Subs:ubd},
where we shall exploit some ideas taken from \cite{BNR}. The case of closed curves will be readily obtained in Sec.~\ref{Subs:closed}.
\smallskip\par We expect that our notion of $p$-rotation $\gk_p(P)$ may be useful from the point of view of numerical analysis, as it allows one to obtain a discretization
of the {$p$-energy} in terms of {\em inscribed polygonals}. In fact, Theorem~\ref{Tmain} may be compared with e.g. \cite{BNR,IB,SSW}, where $\GG$-convergence results are obtained {for}
discrete $p$-curvatures of polygonals converging to the given curve in the {topology} induced by the {\em Fr\'echet distance}, see Definition~\ref{Ddist}.
\par Having in mind possible numerical applications, in Sec.~\ref{Subs:rmks} we shall comment our definition of $p$-rotation of
polygonals from Sec.~\ref{Subs:prot}.
More precisely, we shall see that any reasonable (from a geometric viewpoint) different choice of definition of $p$-rotation of polygonals
(as e.g. by taking the $p$-energy of the optimal piecewise smooth curve satisfying the clamping conditions at the middle points of the segments of $P$)
produces the same relaxed energy $\F_p(\gc)$ we have obtained in our Main Result.
\par In conclusion, our definition seems more fitting from a numerical viewpoint,
since it is well-known that, in general, solutions to Euler's elastica problem with clamping conditions cannot be explicitly computed.\vspace{0.2cm}
\textbf{Funding} The research of D.M. was partially supported by the GNAMPA of INDAM. The research of A.S. was partially supported by the GNSAGA of INDAM.
\section{Physical and geometric-analytical motivations}\label{Sec:phgeo}
In this preliminary section, we briefly discuss the {classical} Euler's elastica problem, giving some examples that are strictly related to our proposal of notion of $p$-rotation of polygonals.
We then report some similar features in the framework of discrete differential geometry,
outlining some possible advantages of our approach to applications through numerical analysis.
\subsection{Euler's elastica}\label{Subs:Euler}
Euler's elastica problem consists in minimizing the total squared curvature,
also known as bending energy, among smooth planar curves $\gc$ of fixed length subject to a clamped (i.e., first order) boundary condition.
The bending energy corresponds to the case $p=2$ of the energy functional
$$ \EE_p(\gc):=\int_\gc \kk^p\,d s\,, \qquad p\geq 1$$
where $s$ denotes the arc-length parameter and $\kk$ the {\em scalar curvature} of the curve $\gc$.
\par In his celebrated study, Euler
derives ODEs for solution curves (i.e., critical points) and moreover classifies the types
of solution curves qualitatively. The solution curves are nowadays called Euler's elastica, see \cite{OEB}.
%
%\subsection{Smooth elastic curves}
%
\par More precisely, let $\gc:\ol I_L\to\RN$ denote the arc-length parameterization $\gc=\gc(s)$ of a rectifiable curve, so that $I_L=(0,L)$ with $L=\calL(\gc)$,
the length of $\gc$.
If $\gc$ is of class $C^2(\ol I_L,\RN)$, the {scalar curvature} at $\gc(s)$ agrees with the norm of the
curvature vector $\gk(s)=\ddot\gc(s)$, whence $\EE_p(\gc)$ agrees with the {\em $p$-energy} functional
\beq\label{TCp} \EE_p(\gc)=\int_0^L\Vert\ddot\gc(s)\Vert^p\,ds \eeq
and when $p=1$ one has $\EE_1(\gc)=\TC(\gc)$, the total curvature.
\par Notice that if the arc-length parameterization of a rectifiable curve $\gc$ belongs to the Sobolev class $W^{2,1}(I_L,\RN)$ or, equivalently,
if the tantrix $\gt(s):=\dot\gc(s)$ belongs to $W^{1,1}(I_L,\SN)$, where $\SN$ is the Gauss hyper-sphere in $\RN$
$$ \SN:=\{y\in\RN\,:\,|y|=1\} $$
then the $p$-energy functional \eqref{TCp} is well-defined, and actually $\EE_p(\gc)<\i$
provided that $\gc\in W^{2,p}(I_L,\RN)$.
In that case, the curve $\gc$ has finite total curvature, $\TC(\gc)<\i$, and finally $\EE_q(\gc)<\i$ for every exponent $1<q<p$.
Moreover, the integral \eqref{TCp} represents the $p$-energy, or the total curvature when $p=1$, also when the tantrix $\gt=\dot\gc$ is continuous and piecewise $C^1$. This corresponds to what we call here the class of {\em piecewise smooth} curves.
%
%\smallskip
\par We report here the result of the computation of the {\em Euler-Lagrange equation} of functionals depending on the scalar curvature of smooth planar curves.
\begin{proposition} For any non-negative smooth function $f$ of the scalar curvature
$\kk$ of smooth curves $\gc$ in $\gR^2$, the Euler equation of the functional
$\F(\gc):=\ds\int_\gc f(\kk)\,ds$ {reads as}
$$ %\beq\label{ELeq}
{\ddot f(\kk)\,\ddot\kk \over |\dot \gc|}+{\dddot
f(\kk)\,\dot\kk^2 \over |\dot \gc|}- {\ddot f(\kk)\,\dot\kk\,(\dot
\gc\bullet \ddot \gc) \over |\dot \gc|^3} + \kk\,\bigl\{\kk\,\dot
f(\kk)-f(\kk)\bigr\}\,|\dot \gc|=0\,. $$%\eeq
\ep
\par By choosing arc-length parameterization, one has $\Vert\dot
\gc\Vert=1$ and $\dot \gc\bullet \ddot \gc=0$, whence the above
equation reduces to the classical one
$$  {\ddot f(\kk)\,\ddot\kk  }+{\dddot f(\kk)\,\dot\kk^2 }+ \kk\,\bigl\{\kk\,\dot f(\kk)-f(\kk)\bigr\}=0 $$
compare \cite[Ch.~1, Sec.~5]{GH}. For $f(\kk)=\e+|\kk|^p$, where
$p>1$ and $\e>0$, so that $\F(\gc)=\e\cdot\calL(\gc)+\EE_p(\gc)$, it takes the simpler form:
\beq\label{ELeq2} p\,|\kk|^{p-2}\,\ddot\kk +p(p-2)\,|\kk|^{p-4}\,\kk\,\dot\kk^2 +
\kk\,|\kk|^p-{\e\,\kk\over (p-1)}=0
\eeq
that when $p=2$ reads as
$$ 2\,\ddot\kk +
\kk^3-\e\,{\kk}=0\,.
$$
\par Now, searching e.g. for smooth closed planar curves with constant
curvature, i.e. for minimal circles of radius $R$, since
$\kk\equiv R^{-1}$ we deduce that equation \eqref{ELeq2} is solved
when $R=R(\e,p)=((p-1)/\e)^{1/p}$.
\par More generally, in presence of first
order boundary conditions, minimizing planar curves in general depend on
the choice of the exponent $p>1$.
\subsection{A minimum problem}\label{Subs:exmin}
For future use, we now consider a minimum problem for the elastic energy under a clamping condition.
Namely, for $\ell>0$ and $\theta\in]0,\pi[$, consider a polygonal $P=P(\ell,\theta)$ in $\gR^2$ given by two segments of length equal to $\ell/2$
and exterior (or turning) angle $\theta$ at the common vertex.
We denote by $\F(\ell,\theta)$ the family of all piecewise smooth planar curves with end points equal to the end points of $P(\ell,\theta)$, with tangent at the end points parallel to the tangents to the polygonal, and with total curvature equal to $\theta$.
If $\gc$ denotes the arc-length parameterization of an element in $\F(\ell,\theta)$, then $\dot\gc$ is continuous and piecewise $C^1$, and we are thus requiring that
$$ \TC(\gc)=\int_0^L \kk(s)\,ds=\theta \,,\qquad L=\calL(\gc)\,,\quad \kk(s)=\Vert\ddot \gc(s)\Vert $$
%where $\kk(s):=\Vert\ddot\gc(s)\Vert$ is the scalar curvature of the curve at $\gc(s)$,
whereas no condition on the length $L$ of the curve is prescribed, as it will be obtained a posteriori.
\par
With the previous notation, one expects that the minimum of the $p$-energy functional \eqref{TCp} in the class $\F(\ell,\theta)$
is attained by a circular arc, independently of the choice of the exponent $p>1$.
\begin{proposition}\label{Pmin1} We have:
$$ \inf \{\EE_p(\gc)\mid \gc\in\F(\ell,\theta)\}=\theta\cdot\Bigl(\frac{\tan(\theta/2)}{\ell/2}\Bigr)^{p-1}\qquad\fa\,p>1 $$
and the infimum is attained by a circular arc {of} radius $R(\ell,\theta):=(\ell/2)\cdot(\tan(\theta/2))^{-1}$. Moreover,
$\ell\,\cos(\theta/2)< \calL(\gc)< \ell$ if $\gc\in\F(\ell,\theta)$.
\ep
\begin{prooof} We first observe that a curve $\gc$ in $\F(\ell,\theta)$ is convex and with support contained in the triangle given by the convex hull of the polygonal $P(\ell,\theta)$. Otherwise, denoting by $\gc$ its arc-length parameterization, {since} the first order boundary conditions read as $\dot\gc(0)\bullet\dot\gc(L)=\arccos\theta$, {we would obtain} $\TC(\gc)>\theta$.

Therefore,
its length $L=\calL(\gc)$ is greater than the distance between the end points of the polygonal $P(\ell,\theta)$, i.e., $L> 2\cos(\theta/2)$,
and lower than the length $\ell$ of the polygonal itself.
Moreover, by {Jensen's} inequality
\beq\label{Je} \Bigl(\frac 1L\,\int_0^L \kk(s)\,ds\Bigr)^p\leq \frac 1L\,\int_0^L \kk(s)^p\,ds \eeq
so that
$$ \int_\gc \kk(s)^p\,d s\geq L^{1-p}\,\theta^p\,. $$
On the other hand, equality holds in \eqref{Je} if and only if $\kk(s)$ is constant, so that $\kk(s)=R(\ell,\theta)^{-1}$, whence $L=R(\ell,\theta)\cdot\theta$ and
$$ { \int_\gc \kk(s)^p\,d s } =\theta\cdot R(\ell,\theta)^{1-p}=\theta\cdot\Bigl(\frac{\tan(\theta/2)}{\ell/2}\Bigr)^{p-1}\,, $$
as required. \epf
\br\label{Ropt} If we modify the polygonal $P$ by requiring that the two segments have different length $\ell_1/2$ and $\ell_2/2$,
where e.g. $\ell_1< \ell_2$, when $p>1$ the energy minimum among piecewise smooth curves satisfying the previous clamping boundary conditions and with total curvature {$\theta$,} is attained by a convex and smooth curve with non-constant curvature.
\par More precisely, one may take among the competitors the circular arc with radius $R(\ell_1/2,\theta)$ attached to a segment of positive length $(\ell_2-\ell_1)/2$.
Now, for any choice of $p>1$ one can always find a convex curve $\gc$ satisfying the clamping boundary conditions and with $p$-energy strictly smaller than the energy given by piecewise smooth curve previously described.
When $\ell_1<\ell_2$, finding the energy minimum is a non-trivial problem, even in the case $p=2$ of the so called bending energy. However, one expects that for small turning angles $\theta$ the energy minimum is comparable to the energy of the circular arc with radius $R(\ell_1/2,\theta)$. \er
\par These formal arguments may justify our choice of {\em $p$-rotation} of polygonal curves $P$, see Figure \ref{fig:pol} on page \pageref{fig:pol} and Definition~\ref{Dprot} below,
that is given by the $p$-energy of the ``optimal" inscribed piecewise smooth curve given by the union of suitable circular arcs.
Moreover, our example in the next section may clarify the situation.
\subsection{Another minimum problem}\label{Subs:exmin2}
Consider now for any $\theta\in]0,\pi[$ a closed polygonal given by a rhombus $P_\theta$ with segments of length one and two opposite external angles equal to $\theta$, and let $\G(\theta)$ denote the class of piecewise smooth and closed curves $\gc$ inscribed in $P_\theta$. We wish to minimize
the $p$-energy functional \eqref{TCp} in the class $\G(\theta)$.
\par Without loss of generality, we may and do assume $\gc$ convex. Moreover, by a symmetry argument it turns out that $\gc$ is tangent to the polygonal $P_\theta$ at four points, whose distance from the nearest between the vertices with external angle $\theta$ is equal to $\l$,
for some $\l\in]0,1[$.
\begin{proposition}\label{Pmin2} For every $p>1$ and $\theta\in]0,\pi[$, the infimum of the problem
$$ \inf \{\EE_p(\gc)\mid \gc\in\G(\theta)\} $$
is attained by a curve $\gc=\gc(\theta,p)$ given by the union of four circular arcs. More precisely, the distance $\l(\theta,p)$ of the tangential points of $\gc(\theta,p)$ from the nearest between the vertices with external angle $\theta$ is equal to
\beq\label{lambda} \l(\theta,p)=\Bigl( 1+\Bigl( \frac{F_p(\pi-\theta)}{F_p(\theta)}\Bigl)^{1/p} \Bigr)^{-1} \eeq
where
$$ F_p(\a):=\a\cdot\tan(\a/2)^{p-1}\,,\quad 0<\a<\pi\,. $$
\ep
\bpf On account of Proposition~\ref{Pmin1}, given the four points in the segments of $P_\theta$ whose distance from the nearest between the vertices with external angle $\theta$ is equal to $\l$, the energy minimizing closed curved $\gc(\l)$ among the ones which are tangential to $P_\theta$ at those fixed four points is given by the union of four circular arcs, two with curvature radius equal to $\l/\tan(\theta/2)$ and total curvature $\theta$, and two with curvature radius equal to $(1-\l)/\tan((\pi-\theta)/2)$ and total curvature $\pi-\theta$. Therefore, one has
$$ \int_{\gc(\l)}\kk^p\,ds=2\,\bigl( {F_p(\theta)}\cdot{\l^{1-p}}+{F_p(\pi-\theta)}\cdot{(1-\l)^{1-p}}\bigr)=:f(\l)\,. $$
\par For each $\theta\in]0,\pi[$ and $p>1$ we have:
$$ f'(\l)=2\,\bigl( {F_p(\theta)}\cdot(1-p)\,\l^{-p}-F_p(\pi-\theta)\cdot(1-p)\,(1-\l)^{-p}\bigr) $$
whence
$$ f'(\l)=0 \iff F_p(\theta)\cdot\l^{-p}=F_p(\pi-\theta)\cdot(1-\l)^{-p} \iff \Bigl(\frac{1-\l}\l\Bigr)^p=\frac{F_p(\pi-\theta)}{F_p(\theta)} $$
which yields to $\l=\l(\theta,p)$ given by \eqref{lambda}, as required.
\epf
\br\label{Ropt2} A part from the case of a square, i.e., when $\theta=\pi/2$, it turns out that the minimal curve is {attained} in correspondence
to tangential points which are not the middle points of the segments of the polygonal, in order to balance the curvature radius of the two couples of inscribed circular arcs.
\par Actually, the energy of the piecewise smooth curve given by four circular arcs tangent to the middle points of the segments of the rhombus $P_\theta$, i.e., the choice of $\l=1/2$ in the proof of Proposition~\ref{Pmin2}, corresponds to our notion of $p$-rotation $\gk_p(P_\theta)$ of the closed polygonal $P_\theta$.
More precisely, on account of Definition~\ref{Dprot} {we get}
$$ \gk_p(P_\theta)= 2^p\,\Bigl( \theta\cdot(\tan(\theta/2))^{p-1}+(\pi-\theta)\cdot(\tan((\pi-\theta)/2))^{p-1} \Bigr) $$
for every $p>1$ and $\theta\in]0,\pi[$, so that when $\theta\neq \pi/2$ {clearly}
$$ \inf \{\EE_p(\gc)\mid \gc\in\G(\theta)\}<\gk_p(P_\theta)\qquad\fa\, p>1\,. $$
\er
\subsection{Discrete elastica}\label{Subs:dg}
In discrete geometry, several definitions are proposed in order to give a discrete analogous to the energy functional $\EE_p(\gc)$ of smooth curves.
Taking for simplicity a closed equilateral polygonal $P$ with $n$ segments equal in length to $\ell$,
and denoting by $\theta_i$ the {\em turning angle} at the  $i$-th vertex $v_i$ of $P$, i.e., the exterior angle between the consecutive segments concurring at $v_i$,
one may take e.g.
\beq\label{Kstar}
\gk_p^*(P):=\sum_{i=1}^n \frac{\theta_i^p}{(\ell/2)^{p-1}}\,,\quad p\geq 1 \eeq
so that for $p=1$ one recovers the {\em rotation} $\gk^*(P)$ of the polygonal, i.e., the sum of the turning angles.
\br\label{Ren} Coming back to Proposition~\ref{Pmin1}, we point out that since $\tan(\theta/2)>\theta/2$, then
$$ \inf \{\EE_p(\gc)\mid \gc\in\F(\ell,\theta)\}>\frac{\theta^p}{\ell^{p-1}}$$
for every $\theta\in]0,\pi[$ and $p>1$, whereas comparing with \eqref{Kstar}
$$ \inf \{\EE_p(\gc)\mid \gc\in\F(\ell,\theta)\}<\frac{\theta^p}{(\ell/2)^{p-1}}$$
provided that $\tan(\theta/2)<\theta$, i.e., for $\theta$ sufficiently small. \er
\par
{In all reasonable definitions of discrete $p$-curvature} available in literature, see e.g. \cite{BNR,CW,IB,MV,SSW}, when $p>1$ the term corresponding to $v_i$ depends on the length of the two consecutive segments concurring at $v_i$, and the corresponding curvature measure is concentrated at the edge points.
\par As we already mentioned in the introduction, in our definition of $p$-rotation from Sec.~\ref{Subs:prot} we view a polygonal as an approximation to a smooth curve and hence spread this curvature out evenly into a smooth density.
\section{Background material and preliminary results}
In this section, we collect some notation and well-known results.
\subsection{Total variation}\label{Subs:tv}
We refer to Secs.~3.1 and 3.2 of \cite{AFP} for the following well-known facts.
\par Let $I=(a,b)\sb\gR$ be a bounded open interval, and $N\in\Nat^+$.
A vector-valued summable function $u:I\to\RN$ is said to be of {\em bounded variation}
if its distributional derivative $Du$ is a finite $\RN$-valued measure in $I$.
\par
The {\em total variation} $|Du|(I)$ of a function $u\in\BV(I,\RN)$
is given by
$$ |Du|(I):=\sup\Bigl\{ \int_I \vf'(s)\,u(s)\,ds \mid \vf\in C^\infty_c(I,\RN)\,,\quad \Vert\vf\Vert_\infty\leq 1 \Bigr\} $$
and hence it does not depend on the choice of the representative in the
equivalence class of the functions that agree $\calL^1$-a.e. in $I$ with $u$, where
$\calL^1$ is the Lebesgue measure.
\par We say that a sequence $\{u_h\}\sb\BV(I,\RN)$ converges to $u\in \BV(I,\RN)$ {\em weakly-$^\ast$ in $\BV$} if $u_h$ converges to $u$ strongly in $L^1(I,\RN)$ and
$\sup_h|Du_h|(I)<\i$. In this case, the lower semicontinuity inequality holds:
$$|Du|(I)\leq\liminf_{h\to\i}|Du_h|(I)\,.$$
If in addition $|Du_h|(I)\to|Du|(I)$, we say that
$\{u_h\}$ {\em strictly converges} to $u$.
\par The {\em weak-$^\ast$ compactness} theorem yields that if $\{u_h\}\sb\BV(I,\RN)$ converges $\calL^1$-a.e. on $I$ to a function $u$,
and if $\sup_h|Du_h|(I)<\i$, then $u\in\BV(I,\RN)$ and a subsequence of $\{u_h\}$ weakly-$^\ast$ converges to $u$.
\par Let $u\in\BV(I,\RN)$. Since each component of $u$ is the difference of two monotone functions, it turns out that $u$ is continuous outside an at most countable set, and that both the left and right limits
$u(s\pm):=\lim_{t\to s^\pm}u(t)$ exist for every $s\in I$. Also, $u$ is an $L^\i$ function that is differentiable $\calL^1$-a.e. on $I$, with derivative $\dot u$ in $L^1(I,\RN)$.
\par The total variation of $u$ agrees with the {\em essential variation} $\Var_{\RN}(u)$, which is equal to the pointwise variation
of any {\em good representative} of $u$ in its equivalence class.
A good (or precise) representative is e.g. given by choosing $u(s)=(u(s+)+u(s-))/2$ at the discontinuity points.
Letting $u_\pm(s):=u(s\pm)$ for every $s\in I$, both the left- and right-continuous functions $u_\pm$ are good representatives.
\par If $u\in\BV(I,\RN)$, the decomposition into the {\em absolutely continuous}, {\em Jump}, and {\em Cantor} parts holds:
$$ Du=D^au+D^Ju+D^Cu\,,\quad |Du|(I)=|D^au|(I)+|D^Ju|(I)+|D^Cu|(I)\,.$$
More precisely, one splits $Du=D^au+D^su$ into the absolutely continuous and singular parts w.r.t. Lebesgue measure $\calL^1$.
The Jump set $J_u$ being the (at most countable) set of discontinuity points of any good representative of $u$, and $\delta_s$ denoting the unit Dirac mass at $s\in I$,
one has:
$$D^a u=\dot u\,\calL^1\,,\qquad D^Ju=\sum_{s\in J_u}[u(s+)-u(s-)]\,\delta_s\,,\qquad D^Cu=D^su\pri(I\sm J_u)\,. $$
\par Also, any $u\in\BV(I,\RN)$ can be represented by $u=u^a+u^J+u^C$, where $u^a$ is a Sobolev function in $W^{1,1}(I,\RN)$, $u^J$ is a Jump function, and $u^C$ is a Cantor function, so that
$$ |D^au|(I)=|Du^a|(I)\,,\qquad |D^Ju|(I)=|Du^J|(I)\,,\qquad |D^Cu|(I)=|Du^C|(I)\,.  $$
As a consequence, $u\in W^{1,1}(I,\RN)$ provided that $D^Ju=0$ and $D^Cu=0$.
\subsection{Length}\label{Subs:l}
Consider a curve $\gc$ in the Euclidean space $\RN$ parameterized by the continuous map
$\gc:\ol I\to\RN$, where $\ol I=[a,b]$.
%, with components $\gc(t)=(\gc^1(t),\ldots,\gc^{N+1}(t))$.
%
Any polygonal curve $P$
{\em inscribed} in $\gc$, say $P \ll \gc$, is obtained by choosing a
finite partition $\D:=\{a=t_0<t_1<\ldots<t_{n-1}<t_{n}=b\}$ of
$\ol I$, say $P=P(\D)$, and letting $P:[a,b]\to\RN$ such that
$P(t_i)=\gc(t_i)$ for $i=0,\ldots,n$, and $P(t)$ affine on
each interval $[t_{i-1},t_i]$. Setting $\gv_i:=\gc(t_i)-\gc(t_{i-1})$ and $\ell_i:=\Vert\gv_i\Vert$, where $\Vert\cdot\Vert$ is the Euclidean norm in $\RN$, the length of the polygonal $P$ is
$\calL(P)=\sum_{i=1}^{n}\ell_i$.
We also denote
$$\mesh \D:=\sup_{1\leq i\leq n}(t_i-t_{i-1})\,,\quad \mesh P:=\sup_{1\leq i\leq n}\ell_i \,. $$
\par The {\em length} $\calL(\gc)$ of the curve $\gc$ is defined by
$$\calL(\gc):=\sup\{\calL(P)\mid P\ll \gc\} $$
and $\gc$ is said to be {\em rectifiable} if $\calL(\gc)<\i$.
\par By uniform continuity, for each $\e>0$ we can find $\d>0$
such that $\mesh P<\e$ if $\mesh \D<\d$ and $P=P(\D)$. As a
consequence, taking $P_h=P(\D_h)$, where
$\{\D_h\}$ is any sequence of partitions of $I$ such that
$\mesh \D_h\to 0$,
we get $\mesh P_h\to 0$ and hence the
convergence $\calL(P_h)\to\calL(\gc)$ of the length functional.
\par The curve $\gc$ is rectifiable if and only if $\gc\in\BV(I,\RN)$, and
in that case
$$\calL(\gc)=\Var_{\RN}(\gc)=|D\gc|(I)\,. $$
Therefore, if $\gc\in C^1([a,b],\RN)$ we get $\calL(\gc)=\int_a^b\Vert\dot\gc(t)\Vert\,dt<\i$.
\bdf\label{Ddist} The {\em Fr\'echet distance} $d(\gc_1,\gc_2)$ between two rectifiable curves is the infimum, over all strictly monotonic
reparameterizations, of the maximum pointwise distance. \edf
Therefore, if $d(\gc_1,\gc_2)=0$, the two curves are equivalent in the following sense: homeomorphic reparameterizations that approach the infimal value zero will
limit to the more general reparameterization that might eliminate or introduce intervals of
constancy, {compare \cite{Su_curv}.}
\par Moreover, if $\{\gc_h\}$ is a sequence of rectifiable curves in $\RN$ such that $d(\gc_h,\gc)\to 0$ as $h\to \i$ for some rectifiable curve $\gc$, then by lower semicontinuity
\beq\label{lsc}\calL(\gc)\leq\liminf_{h\to\i}\calL(\gc_h)\,. \eeq
\subsection{Total curvature}\label{Subs:TC}
We call {\em rotation} $\gk^*(P)$ of a polygonal curve $P$ in $\RN$ the sum of the {\em turning angles} (i.e., the exterior angles) $\theta_i$ between consecutive (and non degenerate)
segments $\gv_i$ and $\gv_{i+1}$.
More precisely, with the previous notation
$$\gk^*(P)=\sum_{i=1}^{n-1} \theta_i\,,\qquad \theta_i:=\arccos\Bigl( \frac{\gv_i\bullet \gv_{i+1}}{\ell_i\cdot\ell_{i+1}}\Bigr) $$
if $\ell_i>0$ for $i=1,\ldots,n$ and the polygonal is open, $\bullet$ denoting the scalar product in $\RN$.
If $P$ is closed, the turning angle between $\gv_n$ and $\gv_1$ is added.
\par
Following Milnor \cite{Mi}, the {\em total curvature} $\TC(\gc)$ of a curve $\gc$ in $\RN$ is defined {by}
$$\TC(\gc):=\sup\{\gk^*(P)\mid P\ll \gc\}\,. $$
Then $\TC(P)=\gk^*(P)$ for each polygonal $P$.
Moreover, if a curve $\gc$ has compact support and finite total curvature, $\TC(\gc)<\i$, then it is a rectifiable curve.
\par Assume now that a rectifiable curve $\gc$ is parameterized by arc-length, so that $\gc=\gc(s)$, with $s\in[0,L]=\ol I_L$, where
$I_L:=(0,L)$ and $L=\calL(\gc)$.
If $\gc$ is smooth and regular, one has $\TC(\gc)=\int_0^L\Vert\gk(s)\Vert\,ds$,
where $\gk(s):=\ddot\gc(s)$ is the curvature vector.
More generally, since $\gc$ is a Lipschitz function, by Rademacher's theorem (cf. \cite[Thm.~2.14]{AFP}) it is differentiable $\calL^1$-a.e. in $I_L$.
Denoting by $\dot f:={d\over ds}f$ the derivative w.r.t.\ arc-length parameter $s$, the tantrix $\gt=\dot\gc$ exists a.e., and actually
$\gt:I_L\to\RN$ is a function of bounded variation. Since moreover $\gt(s)\in\SN$ for a.e. $s$, where $\SN$ is the Gauss hyper-sphere, we shall write $\gt\in\BV(I_L,\SN)$.
The essential variation $\Var_\SN(\gt)$ of $\gt$ in $\SN$ differs from $\Var_{\RN}(\gt)$, as its definition involves the geodesic distance $d_\SN$ in $\SN$ instead of the Euclidean distance in $\RN$.
Therefore, $\Var_{\RN}(\gt)\leq \Var_\SN(\gt)$, and equality holds if and only if $\gt$ has a continuous representative.
More precisely, by decomposing $\gt=\gt^a+\gt^J+\gt^C$, one obtains:
\beq\label{VarSPt}  \Var_\SN(\gt)=\int_0^L|\dot\gt|\,ds+\sum_{s\in J_\gt}d_\SN(\gt(s+),\gt(s-))+|D^C\gt|(I_L)\,. \eeq
Notice also that in the formula for $\Var_{\RN}(\gt)$, that is equal to $|D\gt|(I_L)$, one has to replace in \eqref{VarSPt} the geodesic distance with the Euclidean distance
$\Vert\gt(s+)-\gt(s-)\Vert$ in $\RN$ at each Jump point $s\in J_\gt$.
\par Notice moreover that the Cantor component $D^C\gt$ is non-trivial, in general.
\bex\label{Eangle}
Let e.g. $\gc_u:\ol I\to\gR^2$, where $I=(0,1)$, denote the Cartesian curve $\gc_u(t):=(t,u(t))$ in $\gR^2$
given by the graph of the primitive $u(t):=\int_0^tv(\lambda)\,d\lambda$
of the classical Cantor-Vitali function $v:\ol I\to \gR$
associated to the ``middle thirds" Cantor set. It turns out that $\gt=(1+v^2)^{-1/2}(1,v)$, whence $\gt$ is a Cantor function, i.e., $D^a\gt=D^J\gt=0$, and
$$D\gt(I)=D^C\gt(I)= \int_{I}{1\over(1+v^2)^{3/2}}\,(-v,1)\,d D^Cv\,. $$
Notice that the angle $\omega$ between the unit vectors $(1,0)$ and $\gt$ satisfies $\omega=\arctan v\in\BV(I)$. Therefore, $D\omega(I)=D^C\omega(I)=\int_{I}{1\over 1+v^2}\,d D^Cv $, which yields
$$ |D\omega|(I)=\int_{I}{1\over 1+v^2}\,d |D^Cv|=|D\gt|(I)=\TC(\gc_u)={\p\over 4}\,. $$
\eex
\par The following facts hold:
\ben\item if $P$ and $P'$ are inscribed polygonals to $\gc$ and $P'$ is obtained by adding a vertex in $\gc$ to the vertices of $P$, then $\gk^*(P)\leq\gk^*(P')$\,;
\item if $\gc$ has finite total curvature, for each point $v$ in $\gc$, small open arcs of $\gc$ with an end point equal to $v$ have small total curvature. \een
\par As a consequence, compare \cite{Su_curv}, it turns out that
$\TC(\gc)=\Var_{\SN}(\gt)$, see \eqref{VarSPt}, and that the total curvature of $\gc$ is equal to the limit of $\gk^*(P_h)$ for {\em any} sequence $\{P_h\}$ of polygonals in $\RN$ inscribed in $\gc$ and such that $\mesh P_h\to 0$.
More precisely, if $\gt_h$ is the tantrix of $P_h$, then $\Var_\SN(\gt_h)\to\Var_\SN(\gt)$.
\subsection{Relaxed energies of Cartesian curves}\label{Subs:relcart}
A corresponding relaxed energy functional for Cartesian curves $\gc_u(t)=(t,u(t))$ in $\RN$ was analyzed in \cite{AcMu}, for the ``plastic case" $p=1$,
and in \cite{AcMu2}, for the ``elastic case" $p>1$. Namely, for $u\in C^2(\ol I,\gR^N)$, denote
$$  \EEE_p(u):=\calL(\gc_u)+\EE_p(\gc_u)\,,\qquad \EE_p(\gc_u):=\int_{\gc_u}{\kk_u}^p\, ds  $$
where $\kk_u$ is the scalar curvature of $\gc_u$.
A crucial role is played by the
{\em Gauss map} $\tt_u:I\to\SN$
\beq\label{tau} \tt_u={\dot \gc_u\over \Vert\dot \gc_u\Vert}\,,\qquad \dot
\gc_u=(1,\dot u^1,\ldots,\dot u^N)\,. \eeq
In fact, using that $\Vert\dot \gc_u\Vert\,\kk_u=\Vert\dot\tt_u\Vert$, by the area
formula we get
$$ %\beq\label{TCintro}
\EE_p(c_u)=\int_I\Vert \dot c_u\Vert^{1-p}\,\Vert\dot
\tt_u\Vert^p\,dt\,,\quad p\geq 1\,. $$
%\eeq
%
\par
For any summable function $u\in L^1(I,\gR^N)$, the relaxed energy is defined by
$$ \ol\EEE_p(u):=\inf\{\liminf_{h\to\i}
\EEE_p(u_h)\mid \{u_h\}\sb C^2(I,\gR^N)\,,\,\, u_h\to
u\,\,{\text{in}}\,\, L^1(I,\gR^N)\} $$
and clearly if $\ol\EEE_p(u)<\i$ for some $p>1$, then $\ol\EEE_1(u)<\i$.
\par Let now $u$ be a continuous functions $u\in C^0(\ol I,\gR^N)$, so that $\gc_u$ is a compactly supported Cartesian curve.
\par Condition $\ol\EEE_1(u)<\i$ yields that $\gc_u$ is rectifiable, whence the Gauss map $\tt_u$ is well-defined $\calL^1$-a.e. in $I$ by \eqref{tau},
but in terms of the approximate gradient of $u$.
In addition, $\tt_u$\, is a function of bounded variation in $\BV(I,\SN)$, the total variation of $\tt_u$ agrees with the total curvature of the Cartesian curve, and
$$  \ol\EEE_1(u)=\calL(\gc_u)+\TC(\gc_u)\,,\quad \TC(\gc_u)=|D\tt_u|(I)\,. $$
\par If $\ol\EEE_p(u)<\i$ for some exponent $p>1$, then $\tt_u$ is a {\em special function of bounded variation}, i.e., $D^C\tt_u=0$.
Therefore, the planar curve from Example~\ref{Eangle} satisfies $\ol\EEE_p(u)=\i$ for each $p>1$, whereas $\calL(\gc)=\int_0^1\sqrt{1+v^2}\,dt$ and with
$\gt=\tt_u$ we get $\TC(\gc_u)=\pi/4$.
\par
Moreover, in codimension $N=1$, {it turns out that} a Cartesian curve $\gc_u$ with finite relaxed elastic energy, i.e., $\ol\EEE_p(u)<\i$ for some $p>1$, cannot have
corner points, compare \cite{DFLM}. In fact, $D^J\tt_u=0$ and hence $\tt_u\in W^{1,1}(I,\Su)$.
Most importantly, if $\ol\EEE_p(u)<\i$ the arc-length parameterization $\gc:\ol I_L\to\gR^2$ of the curve $\gc_u$ is a Sobolev map in $W^{2,p}(I_L,\gR^2)$, and actually
\beq\label{Erel}  \ol\EEE_p(u)=\int_0^L\bigl(1+\Vert\ddot\gc(s)\Vert^p\bigr)\,ds<\i\,. \eeq
In high codimension $N\geq 2$, corner points may appear. However, roughly speaking, the relaxation process generates a rectifiable curve $\gc$ for which formula \eqref{Erel} continues to hold.
\section{p-curvature of non-smooth curves}\label{Sec:pcurv}
In this section, we introduce our notion of {{\em $p$-rotation} of polygonals and of {\em $p$-curvature}} of non-smooth curves, outlining its main properties.
\subsection{p-rotation of a polygonal}\label{Subs:prot}
Let $P$ be an open polygonal in $\RN$ given by $n$ non-degenerate and consecutive segments $\gv_i$ of length $\calL(\gv_i)=\ell_i>0$, for $i=1,\ldots,n$.
Assume that the turning angle $\t_i$ between the consecutive segments $\gv_{i}$ and $\gv_{i+1}$ concurring at the vertex $v_i$ is lower than $\pi$, for
$i=1,\ldots,n-1$.
\par Let $r_i:=\min\{\ell_i,\,\ell_{i+1}\}$ and denote by $t_{i}^-$ and $t_i^+$ the points in the segments $\gv_i$ and $\gv_{i+1}$, respectively,
whose distance to $v_i$ is equal to $r_i/2$.
\par
Also, if $\theta_i>0$, we let $\g_i$ denote the {\em oriented circular arc} with initial point $t_{i}^-$, final point $t_i^+$, and with tangent parallel to $\gv_i$ and
$\gv_{i+1}$ at $t_{i}^-$ and $t_i^+$, respectively.
When $\theta_i=0$, then $\g_i$ is the segment between $t_{i}^-$ and $t_i^+$, a degenerate ``circular arc" with zero curvature.

See Figure \ref{fig:pol} on page \pageref{fig:pol} for reference on this construction.
\bdf\label{DgP} With the previous notation, we denote by $\g(P)$ the piecewise smooth curve that parameterizes consecutively the arc $\g_1$, the segment between $t_1^+$ and $t_2^-$, the arc $\g_2$, the segment between $t_2^+$ and $t_3^-$, and so on until we get to the final arc $\g_{n-1}$.
\edf
\br\label{Rclosed}
If the polygonal $P$ is closed, the definition is modified in a straightforward way by also considering the angle at the end points $v_n=v_0$, so that $\g(P)$ becomes a piecewise smooth and closed curve.
\er
\par It is readily seen that the curve $\g(P)$ is rectifiable and with length $\calL(\g(P))$ lower than the length of $P$. Moreover, the arc-length parameterization $\gc(P):[0,\calL(\g(P))]\to\RN$ of $\g(P)$ is a piecewise smooth function. The scalar curvature $\kk_{\gc(P)}$ of $\gc(P)$ is equal to zero in correspondence to the segments, and equal to the reciprocal of the radius of the circle completing $\g_i$, at the points inside $\g_i$.
\bdf\label{Dprot} With the previous notation, for every $p\geq 1$ we call {\em $p$-rotation} of $P$ the number
$$ \gk_p(P):=\int_{\gc(P)}{\kk_{\gc(P)}}^p\,ds  =\int_0^{\calL(\g(P))} (\kk_{\gc(P)}(s))^p\,ds\,. $$
\edf
\par When $p=1$, it is readily checked that $\gk_1(P)$ is equal to the rotation $\gk^*(P)$ of the polygonal.
In addition, if a turning angle of $P$ is equal to $\p$, we let $\gk_p(P):=+\infty$ for each $p>1$.
\bex If $P$ is an equilateral closed polygonal $P$ with $n$ segments of length equal to $\ell$, denoting by $\theta_i$ the turning angle at the $i$-th vertex $v_i$, we obtain
$$\gk_p(P)=\sum_{i=1}^n \theta_i\cdot\Bigl( \frac{\tan(\theta_i/2)}{\ell/2}\Bigr)^{p-1}=\frac 1{(\ell/2)^{p-1}}\sum_{i=1}^n \theta_i\cdot
({\tan(\theta_i/2)})^{p-1} $$
and hence, comparing with \eqref{Kstar}, since $\theta_i>0$ for some $i$, in general we get
$$ \gk_p^*(P)< \gk_p(P)\,. $$
\par In particular, if $P=P_{n,\ell}$ is a regular $n$-agon in the plane with edges of length $\ell$, then $\g(P_{n,\ell})$ is its inscribed circle,
and one gets
$$ \gk_p(P_{n,\ell})= 2\pi\Bigl( (2/\ell)\cdot\tan(\pi/n)\Bigr)^{p-1}\qquad p\geq 1\,.$$
Notice that if the circle $\g(P_{n,\ell})$ is assumed to have perimeter equal to one, then $\ell=\ell_n=\pi^{-1}\tan(\pi/n)$ and for each $n\geq 3$ we get
$\gk_p(P_{n,\ell_n})=(2\pi)^p\to +\i$ if $p\to +\i$. \eex
\smallskip\par Now, if $P$ and $P'$ are two polygonals inscribed to a curve $\gc$, and $P'$ is obtained by adding a vertex in $\gc$ to the vertices of $P$, differently from the case $p=1$, where $\gk_1(P)$ is the rotation $\gk^*(P)$ of the polygonal,
when $p>1$ in general one cannot compare the $p$-rotation $\gk_p(P)$ of $P$ with the $p$-rotation $\gk_p(P')$ of $P'$.
\bex\label{Emonot} If $P$ is obtained by four collinear and consecutive vertices $v_i$, $i=1,\ldots,4$, and $P'$ by adding a fifth non-collinear vertex $w$ between $v_2$ and $v_3$, clearly $0=\gk_p(P)<\gk_p(P')$ for all $p\geq 1$.
Instead, if $P$ is the edge of a square with side length equal to two, then $\gk_p(P)=2\pi$ for every $p\geq 1$. Adding an external vertex $w$ to the square in such a way that the triangle with vertices $v_1,w,v_2$ is equilateral (where $v_1$ and $v_2$ are the nearest vertices of $P$ to $w$) we obtain a polygonal $P'$ such that $\gk_p(P)=3\pi/2+\pi\,(\tan(\pi/8))^{p-1}/2$, whence $\gk_p(P)>\gk_p(P')$ for every $p>1$.
\eex
\par The lack of validity of a monotonicity formula justifies our approach in the definition of a $p$-curvature functional.
\subsection{The p-curvature functional}\label{Subs:pcurv}
Due to the drawback outlined by Example~\ref{Emonot}, differently to the case $p=1$, one cannot introduce a notion of $p$-curvature for non-smooth curves in terms of supremum of the $p$-rotation of inscribed polygonals.
We thus follow the classical approach that goes back to Lebesgue-Serrin definition of relaxed functional, and make use of the notion of {\em modulus} by Alexandrov-Reshetnyak \cite{AR}.
\par
In the sequel, we always assume that $\gc$ is a rectifiable curve in $\RN$
%with finite total curvature, $\TC(\gc)<\i$,
parameterized in arc-length.
The modulus $\m_\gc(P)$ of a polygonal $P$ inscribed in $\gc$, say $P\ll\gc$, is the maximum of the diameter of the arcs of $\gc$ determined by the consecutive vertices in $P$.
%In the same spirit as Lebesgue-Serrin relaxed functional, we introduce the following:
%
\bdf\label{DEp} We call {\em $p$-curvature} $\F_p(\gc)$ of a rectifiable curve $\gc$ in $\RN$ the number
$$%\beq\label{DTAT}
 \F_p(\gc):=\inf\Bigl\{\liminf_{h\to\i}\gk_p(P_h)\mid \{P_h\}\ll\gc\,,\,\,\m_\gc(P_h)\to 0\Bigr\} $$
% \F_p(\gc):=\lim_{\e\to 0^+}\sup\{ \gk_p(P) \mid P\ll \gc \,,\,\,\m_\gc(P)<\e\}
%\eeq
where $\gk_p(P)$ is the $p$-rotation of $P$, see Definition~\ref{Dprot}.
\edf
\par Notice that if $p=1$, we get
\beq\label{F1TC} \F_1(\gc)=\TC(\gc) \eeq
so that if $\gc$ is a polygonal curve itself, we get $\F_1(P)=\gk^*(P)<\i$.
\par However, for a polygonal curve $P$ with positive rotation, $\gk^*(P)>0$, we clearly have
$$ \F_p(P)=+\infty\qquad\fa\,p>1\,. $$
In fact, if $P$ has a positive turning angle {$\t\in]0,\pi[$} at a vertex $v$, for each $\e>0$ sufficiently small, we can define an inscribed polygonal $P_\e$ containing the consecutive vertices $u,v,w$, where $u$ and $w$ lie at a distance $2\e$ from $v$ on the edges concurring at $p$.
The contribution to the $p$-rotation of $P_\e$ near the vertex $v$ is comparable to $\e^{1-p}\cdot\theta\,(\tan(\t/2))^{p-1}$, whence it diverges to $+\i$ as $\e\to 0^+$.
\par With a similar computation, it turns out that $\F_p(\gc)=+\infty$ for every $p>1$, if $\gc$ is any piecewise smooth curve with at least one corner point. This is coherent with the physical interpretation: an elastic rod
needs infinite bending energy in order to produce a corner.
%
%\par Notice also that if $\F_p(\gc)<\i$,
%for any sequence $\{P_h\}$ of polygonal curves inscribed in $\gc$ and satisfying $\m_\gc(P_h)\to 0$, one has $\sup_h\gk_p(P_h)<\i$,
%one can also find an optimal sequence as above in such a way that $\gk_p(P_h)\to\F_p(\gc)$.
%

Moreover, we have:
\begin{proposition}\label{Ppq} Let $\gc$ be a rectifiable curve in $\RN$ such that $\F_p(\gc)<\i$ for some $p>1$. Then $\F_q(\gc)<\i$ for all $1\leq q<p$, and
$$ \F_q(\gc)\leq\calL(\gc)+\F_p(\gc)\,. $$
In particular, $\gc$ has finite total curvature, $\TC(\gc)<\i$.
\ep
\bpf If $1\leq q< p$, then $t^q\leq 1+t^p$ for every $t>0$. For every polygonal curve $P$ inscribed in $\gc$, according to Definition~\ref{Dprot}, and recalling that $\calL(\g(P))\leq\calL(P)$ we thus get
$$ \gk_q(P)\leq\calL(P)+\gk_p(P) $$
where $\gk_1(P)=\gk^*(P)$. Since $\calL(P)\leq\calL(\gc)$, the claim follows on account of Definition~\ref{DEp}.
\epf
\section{Main results}
Let $\gc$ be a rectifiable and open curve in $\RN$, and let $\gc:\ol I_L\to\RN$ denote its arc-length parameterization, so that $I_L=(0,L)$ and $L=\calL(\gc)$. We have noticed that for $p=1$, our {$p$-curvature functional} from Definition~\ref{DEp} agrees with the total curvature, see \eqref{F1TC}. On the other hand, it is well-known that $\gc$ has finite total curvature if and only if the tantrix $\gt=\dot\gc$ is a function of bounded variation,
$\gt\in\BV(I_L,\SN)$, and in that case $\TC(\gc)=\Var_\SN(\gc)$, compare \eqref{VarSPt}. Therefore, if in particular $\gt$ is a Sobolev map in $W^{1,1}(I_L,\SN)$ we get $\TC(\gc)=\int_0^L\Vert\ddot\gc(s)\Vert\,ds$.
\par A completely different situation occurs when $p>1$.
In accordance with the phenomena observed in the relaxation process for Cartesian curves, we shall prove that a rectifiable and open curve has finite $p$-curvature for some $p>1$ if and only if its tantrix $\gt$ is a Sobolev map in $W^{1,p}(I_L,\SN)$.
More precisely, we have:
\bt\label{Topen} Let $\gc$ be a rectifiable and open curve in $\RN$ parameterized in arc-length. Then for every exponent $p>1$
$$ \F_p(\gc)<\i \iff \gc\in W^{2,p}(I_L,\RN) $$
see Definition~$\ref{DEp}$, and in that case
$$ \F_p(\gc)=\int_0^{L}\Vert\ddot \gc(s)\Vert^p\,ds\,. $$ \et
\par Our Main Result is a direct consequence of Theorems~\ref{Tlb} and \ref{Tub} below. In the first one, we obtain the energy lower bound
and the membership of the tantrix $\gt$ to the Sobolev class $W^{1,p}(I_L,\SN)$. In the second one, we obtain the energy upper bound.
We then consider the case of closed curves, see Corollary~\ref{Cclosed}.
\subsection{Energy lower bound and Sobolev regularity}\label{Subs:lbd}
In this section, we prove the following:
\bt\label{Tlb} Let $\gc$ be a rectifiable and open curve in $\RN$ parameterized in arc-length such that $\F_p(\gc)<\i$ for some $p>1$.
Then $\gc\in W^{2,p}(I_L,\RN)$ and
\beq\label{int3} \int_0^{L} \Vert \ddot\gc(s)\Vert^p\,ds\leq \F_p(\gc)<\i\,. \eeq
%$$ \F_p(\gc)=\int_0^{I_L}\Vert\ddot \gc(s)\Vert^p\,ds\,. $$
\et
\bpf Let $\{P_h\}$ denote an optimal sequence of polygonal curves inscribed in $\gc$, i.e., satisfying $\m_\gc(P_h)\to 0$ and $\gk_p(P_h)\to\F_p(\gc)$. For each $h$, let $\gc(P_h):[0,L_h]\to\RN$ denote the arc-length parameterization of the curve $\g(P_h)$ from Definition~\ref{Dprot}, where $L_h:=\calL(\g(P_h))$, and let $\g_h:[0,L]\to\RN$ given by $\g_h(s):=\gc(P_h)(L_hs/L)$, where $L:=\calL(\gc)$.
By piecewise smoothness, a part from a finite set of points one has $\kk_{\gc(P_h)}(\l)=\Vert\ddot \gc(P_h)(\l)\Vert$ for $\l\in[0,L_h]$ and
$\ddot\g_h(s)=(L_h/L)^2\cdot\ddot\gc(P_h)(\l)$ for $s\in [0,L]$, with $\l=L_hs/L$. Therefore,
\beq\label{int1} \gk_p(P_h)=\int_0^{\calL(\g(P))} (\kk_{\gc(P)}(\l))^p\,d\l=\Bigl( \frac{L}{L_h} \Bigr)^{2p-1}\int_0^{L} \Vert \ddot\g_h(s)\Vert^p\,ds\,. \eeq
\par Now, we have $d(\g(P_h),P_h)\leq\m_c(P_h)$ for every $h$, whereas $d(P_h,\gc)\to 0$, see Definition~\ref{Ddist}. Since $\m_c(P_h)\to 0$, we get $d(\g(P_h),\gc)\to 0$, whence by \eqref{lsc} we infer that
$$ \calL(\gc)\leq\liminf_{h\to\i}\calL(\g(P_h))\,. $$
Using that $\calL(\g(P_h))\leq\calL(P_h)\leq\calL(\gc)$ for every $h$, we deduce that $L_h\to L$.
\par As a consequence, recalling that $\gk_p(P_h)\to\F_p(\gc)$, by \eqref{int1} we obtain:
\beq\label{int2} \lim_{h\to\i}\int_0^{L} \Vert \ddot\g_h(s)\Vert^p\,ds=\F_p(\gc)\,. \eeq
\par Since $p>1$, the sequence $\{\dot\g_h\}$ strongly converges in $W^{1,1}$ to some function $v\in W^{1,1}(I_L,\RN)$. Using that $\g_h$ converges to the Lipschitz function $\gc$ strongly in $L^1(I_L,\RN)$, we get $v=\dot\gc$ a.e., whence possibly passing to a (not relabeled) subsequence, $\{\dot\g_h\}$ converges to $\dot\gc$ weakly in $W^{1,p}(I_L,\RN)$. In particular, $\dot \gc\in W^{1,p}(I_L,\SN)$ and by lower semicontinuity
$$ \int_0^{L} \Vert \ddot\gc(s)\Vert^p\,ds\leq \liminf_{h\to\i}\int_0^{L} \Vert \ddot\g_h(s)\Vert^p\,ds $$
so that by \eqref{int2} we get \eqref{int3}, as required. \epf
%
%Therefore, since $\TC(\g_h)=\TC(P_h)$, where $\TC(P_h)\to\TC(\gc)$ and
%$$ \TC(\g_h)=\frac{L}{L_h}\int_0^{L} \Vert \ddot\g_h(s)\Vert\,ds\,,\quad \TC(\gc)=\int_0^{L} \Vert \ddot\gc(s)\Vert\,ds $$
%we infer that $\{\dot\g_h\}$ strictly converge to $\dot\gc$ in $\BV$, i.e., $\dot\g_h\to \dot\gc$ in $L^1(I_L,\RN)$ and
%$$ \lim_{h\to\i}\int_0^{L} \Vert \ddot\g_h(s)\Vert\,ds= \int_0^{L} \Vert \ddot\gc(s)\Vert\,ds\,. $$
%
%\br Following the previous proof, in a similar way we obtain that for any sequence $\{P_h\}$ of polygonal curves inscribed in $\gc$ such that $\m_\gc(P_h)\to 0$ and $\sup_h\gk_p(P_h)<\i$, then
%$$\int_0^L\Vert\ddot\gc(s)\Vert^p\,dx\leq\liminf_{h\to\i}\gk_p(P_h)\,. $$
%
%Therefore, we obtain:
%
%\bc If $\gc$ is rectifiable curve in $\RN$ parameterized in arc-length and such that $\F_p(\gc)<\i$ for some $p>1$, then $\gc\in W^{2,p}(I_L,\RN)$ and
%
%$$
% \int_0^{L} \Vert \ddot\gc(s)\Vert^p\,ds\leq \inf\Bigl\{\liminf_{h\to\i}\gk_p(P_h)\mid \{P_h\}\,,\,\,\m_\gc(P_h)\to 0\Bigr\}<\i\,. $$
%
%\ec \er
%
\subsection{Energy upper bound}\label{Subs:ubd}
Using some ideas taken from \cite{BNR}, we now obtain the energy upper bound.
\bt\label{Tub} Let $\gc$ be a rectifiable and open curve in $\RN$, parameterized in arc-length. If $\gc$ belongs to $W^{2,p}(I_L,\RN)$ for some $p>1$,
there exists a sequence $\{P_h\}$ of polygonal curves inscribed in $\gc$ such that $\m_\gc(P_h)\to 0$ and
\beq\label{ub} \liminf_{h\to\i}\gk_p(P_h)\leq \int_0^{L} \Vert \ddot\gc(s)\Vert^p\,ds<\i\,. \eeq
In particular, by Definition~$\ref{DEp}$ we have:
$$ \F_p(\gc)\leq \int_0^{L} \Vert \ddot\gc(s)\Vert^p\,ds\,.$$
\et
\bpf We divide the proof in three steps. Firstly, we choose the inscribed polygonals $P_\e$, for $\e>0$ small. Secondly, we make use of some estimates from \cite{BNR} in order to obtain a lower bound of the integral $\int_0^{L} \Vert \ddot\gc(s)\Vert^p\,ds$ in terms of the {$p$-rotation $\gk_p(P_\e)$.}
Finally, we prove the angle estimate from Lemma~\ref{Langle} below.
\smallskip\par\noindent{\sc Step 1:} Let $\gc:[0,L]\to\RN$ be the arc-length parameterization of the curve.
%Choose $0<\e<1$ small.
By absolute continuity, for each $\y>0$ small there exists $\d=\d(\y)>0$ such that if $I\sb[0,L]$ is an interval with length $|I|<\d$, then
$\int_I\Vert\ddot\gc\Vert\,dt<\y$.
%2\,\arccos(1/(1+\e))$.
%
Therefore, if $0<\a<\be< L$ and $\be-\a<\d$, it turns out that $\Vert\dot\gc(\be)-\dot\gc(\a)\Vert\leq \y$ and hence the angle $\theta$ between the unit vectors
$\dot\gc(\a)$ and $\dot\gc(\be)$ is smaller than $2\,\arcsin(\y/2)$.
\par We also notice that if $\gc(\be)-\gc(\a)=\ell$, then any curve with end points $\gc(\a)$ and $\gc(\be)$ and with total curvature $\theta$ has
length $\wid\ell$ lower than $\ell\cdot(\cos(\theta/2))^{-1}$. In particular, we have:
$$ \theta\leq 2\,\arcsin(\y/2) \Longrightarrow \wid\ell\leq \frac {\y}{\sqrt{1-(\y/2)^2}}\,. $$
\par We now fix $0<\e<1$ small and in correspondence we choose $\y_1=\y_1(\e)>0$ so that
$$ \frac {\y_1}{\sqrt{1-(\y_1/2)^2}}\leq \e\,. $$
In addition, we choose $\y_2=\y_2(\e)>0$ small in such a way that
$$ 0\leq 2\theta\leq 2\arcsin(\y_2/2) \Longrightarrow \left\{ \ba{l}2\,\tan(\theta /2)\leq (1+\e)\cdot\theta \\ 2(1-\cos\theta)\geq (1-\e)\cdot{\theta^2}\,. \ea \right. $$
Finally, we let $\y(\e):=\min\{\y_1(\e),\,\y_2(\e)\}$ and define $\d_\e=\d(\y(\e))$ as above.
\par We now may and do choose the greatest number $\ell$ with $\ell\leq \d_\e$ such that we can find an equilateral poligonal $P_\e$ inscribed in $\gc$ and with edge length equal to $\ell$.

More precisely, there exists $n\in\Nat^+$ and $0=t_0<t_1<\cdots<t_{n-1}<t_n=L$ such that
$$ t_i=\min\{t\in[t_{i-1},L]\,:\, \Vert\gc(t)-\gc(t_{i-1})\Vert = \ell\}\qquad\fa\, i=1,\ldots,n $$
and also, by the previous construction,
\beq\label{ell} \ell\leq (t_i-t_{i-1})\leq \ell\cdot(1+\e) \qquad\fa\,i=1,\ldots,n\,. \eeq
\par We then denote by $P_\e$ the equilateral polygonal inscribed in $\gc$ and with consecutive vertices $\gc(t_i)$, for $i=0,\ldots,n$.
Letting $\gv_i:=\gc(t_i)-\gc(t_{i-1})$, we have $\Vert\gv_i\Vert=\ell$ for $i=1,\ldots,n$, and $\calL(P_\e)=n\,\ell$.
Moreover, by Definition~\ref{Dprot} the $p$-rotation of the equilateral polygonal $P_\e$ is
$$ \gK_p(P_\e)={(\ell/2)^{1-p}}\sum_{i=1}^{n-1}\theta_i\cdot\Bigl(\tan\frac{\theta_i}2\Bigr)^{p-1} $$
where $\theta_i$ is the turning angle of $P_\e$ at the vertex {$\gc(t_i)$}.

\par Since by our construction
$$ \tan\frac{\theta_i}2\leq(1+\e)\,\frac{\theta_i}2\qquad\fa\,i=1,\ldots,n-1 $$
we can estimate the $p$-rotation of $P_\e$ as follows:
\beq\label{est2} \gk_p(P_\e)\leq(1+\e)^{p-1}\sum_{i=1}^{n-1} \frac{\theta_i^p}{\ell^{p-1}}\,. \eeq
\par Finally, we have also obtained the uniform lower bound:
\beq\label{est3} 2(1-\cos\theta_i)\geq (1-\e)\,{\theta_i}^2 \qquad\fa\,i=1,\ldots,n-1\,. \eeq
{\sc Step 2:} Let now $P_\e(s):[0,n\,\ell]\to\RN$ denote the arc-length parameterization of $P_\e$, so that if $s_i=\ell\cdot i$, then
$P_\e(s_i)=\gc(t_i)$ for $i=0,\ldots,n$.
\par Following the lines of the proof of Lemma 7 in Appendix 2 of \cite{BNR}, we consider the piecewise linear homeomorphism {$\psi:[0,n\,\ell]\to[0,L]$}
such that {$\psi(s_i)=t_i$} for $i=0,\ldots,n$ and $\psi_{\vert [s_{i-1},s_i]}$ is affine for $i\geq 1$.
\par One clearly has
$$ \int_0^L\Vert\ddot\gc(t)\Vert^p\,dt \geq \sum_{i=1}^{n-1}\frac 1\ell\,\int_0^\ell\Bigl( \int_{\psi(s_{i-1}+a)}^{\psi(s_{i}+a)} \Vert\ddot\gc(t)\Vert^p\,dt \Bigr)\,da $$
where for $i=1,\ldots,n-1$ by Jensen's inequality one obtains the estimate
$$ \frac 1\ell\,\int_0^\ell\Bigl( \int_{\psi(s_{i-1}+a)}^{\psi(s_{i}+a)} \Vert\ddot\gc(t)\Vert^p\,dt \Bigr)\,da \geq
\bigl( \max\{t_{i+1}-t_i,t_i-t_{i-1}\}\bigr)^{1-p}\cdot \Bigl\Vert \frac{\gv_{i+1}}{t_{i+1}-t_i}- \frac{\gv_{i}}{t_i-t_{i-1}} \Bigr\Vert^p\,. $$
Assuming e.g. $t_{i+1}-t_i\geq t_i-t_{i-1}$, one gets
$$ \frac 1\ell\,\int_0^\ell\Bigl( \int_{\psi(s_{i-1}+a)}^{\psi(s_{i}+a)} \Vert\ddot\gc(t)\Vert^p\,dt \Bigr)\,da \geq
\frac{\ell^p}{(t_{i+1}-t_i)^{2p-1}}\cdot \Bigl\Vert \frac{\gv_{i+1}}\ell-\frac{\gv_{i}}\ell +\frac{\gv_{i}}\ell\Bigl( 1-\frac {t_{i+1}-t_i}{t_i-t_{i-1}}\Bigr) \Bigr\Vert^p\,. $$
\par Now, by using \eqref{ell} we obtain the lower bound for the first term:
$$ \frac{\ell^p}{(t_{i+1}-t_i)^{2p-1}}\geq \frac{\ell^{1-p}}{(1+\e)^{2p-1}}\,.$$
\par Moreover, the following angle estimate for the second term holds true:
\bl\label{Langle}
With the previous notation, for $i=1,\ldots,n-1$
$$ \Bigl\Vert \frac{\gv_{i+1}}\ell-\frac{\gv_{i}}\ell +\frac{\gv_{i}}\ell\Bigl( 1-\frac {t_{i+1}-t_i}{t_i-t_{i-1}}\Bigr) \Bigr\Vert^2 \geq
2(1-\cos\theta_i)$$
where, we recall, $\theta_i$ is the turning angle of $P_\e$ at the vertex $\gc(t_i)$. \el
Since the proof of a similar estimate in \cite{BNR} is omitted, we will demonstrate Lemma~\ref{Langle} in Step 3 below.
\par By Lemma~\ref{Langle} and the uniform lower bound \eqref{est3}, we obtain
$$ \Bigl\Vert \frac{\gv_{i+1}}\ell-\frac{\gv_{i}}\ell +\frac{\gv_{i}}\ell\Bigl( 1-\frac {t_{i+1}-t_i}{t_i-t_{i-1}}\Bigr) \Bigr\Vert^p \geq
(1-\e)^{p/2}\cdot{\theta_i^p} \quad\fa\,i=1,\ldots,n-1\,.$$
\par As a consequence, putting the terms together we find
$$\int_0^L\Vert\ddot\gc(t)\Vert^p\,dt \geq \sum_{i=1}^{n-1} \frac{\ell^{1-p}}{(1+\e)^{2p-1}}\cdot (1-\e)^{p/2}\cdot{\theta_i^p} $$
so that by \eqref{est2} we get
$$ \gk_p(P_\e)\leq \frac{(1+\e)^{3p-2}}{(1-\e)^{p/2}}\cdot\int_0^L\Vert\ddot\gc(t)\Vert^p\,dt  $$
and hence the assertion readily follows by letting $P_h:=P_{\e_h}$ for a suitable decreasing sequence $\e_h\searrow 0$.
\adl\par\noindent{\sc Step 3:} It remains to prove the angle estimate in Lemma~\ref{Langle}.
\par By \eqref{ell}, and recalling that we assumed $t_{i+1}-t_i\geq t_i-t_{i-1}$, we have:
$$ \Bigl( 1-\frac {t_{i+1}-t_i}{t_i-t_{i-1}}\Bigr)=-\s $$
for some $0\leq\s\leq \e$, whence we can write
$$ d:=\Bigl\Vert \frac{\gv_{i+1}}\ell-\frac{\gv_{i}}\ell +\frac{\gv_{i}}\ell\Bigl( 1-\frac {t_{i+1}-t_i}{t_i-t_{i-1}}\Bigr) \Bigr\Vert =
\Bigl\Vert \frac{\gv_{i+1}}\ell-\frac{\gv_{i}}\ell \cdot(1+\s) \Bigr\Vert \,. $$
\par We now apply Carnot's theorem to the triangle with sides of length $d$, $L_+$ and $L_-$, where
$$ L_+:=\Bigl\Vert \frac{\gv_{i+1}}\ell\Bigr\Vert\,,\qquad L_-:=\Bigl\Vert \frac{\gv_{i}}\ell \cdot(1+\s) \Bigr\Vert $$
so that the opposite angle to the side of length $d$ is equal to $\theta_i$, obtaining:
$$ d^2=L_+^2+L_-^2-2\,L_+L_-\,\cos\theta_i\,. $$
%=2\,\Bigl( \frac{L_+^1+L_-^2}2-L_+L_-\,\cos\theta_i\Bigr)\,.
%
\par By the estimate \eqref{ell}, and recalling that $0\leq\s\leq\e<1$, we have:
$$ 1\leq L_+\leq (1+\e)\,,\qquad 1+\s\leq L_-\leq(1+\e)\cdot(1+\s)<(1+3\e) $$
so that
$$ L_+=1+\s_+\,,\qquad L_-=1+\s_- $$
with $0\leq\s_+\leq\e$ and $0\leq\s_-\leq 3\e $, whence we re-write
$$\ba{rl} d^2= & (1+\s_+)^2+(1+\s_-)^2-2\,(1+\s_+)(1+\s_-)\,\cos\theta_i \\
= & 2(1+\s_++\s_-)+(\s_+^2+\s_-^2)-2\,(1+\s_++\s_-+\s_+\s_-)\,\cos\theta_i \\
= & 2(1+\s_++\s_-)\,(1-\cos\theta_i)+\bigl( \s_+^2+\s_-^2-2\,\s_+\s_-\,\cos\theta_i \bigr)\,.
 \ea$$
Since $\s_+\geq 0$, $\s_-\geq 0$, and
$$ \bigl( \s_+^2+\s_-^2-2\,\s_+\s_-\,\cos\theta_i \bigr)\geq (\s_+-\s_-)^2\geq 0 $$
we get
$$ d^2\geq 2\,(1-\cos\theta_i) $$
and the proof is complete. \epf
\subsection{The case of closed curves}\label{Subs:closed}
For closed curves, we readily obtain the following
\bc\label{Cclosed} Let $\gc$ be a rectifiable and closed curve in $\RN$ parameterized in arc-length. Then for every exponent $p>1$
$$ \F_p(\gc)<\i \iff \gc\in W^{2,p}(I_L,\RN)\quad{\text{and}}\quad \dot\gc(0)=\dot\gc(L) $$
and in that case
$$ \F_p(\gc)=\int_0^{L}\Vert\ddot \gc(s)\Vert^p\,ds\,. $$ \ec
\bpf Coming back to the proof of Theorem~\ref{Tlb}, for closed rectifiable curves $\gc$, so that $\gc(0)=\gc(L)$, the definition of $p$-rotation of inscribed polygonals is modified on account of Remark~\ref{Rclosed}. Therefore, by condition $\F_p(\gc)<\i$ we obtain again the energy lower bound \eqref{int3}, the membership of $\gc$ to the Sobolev space $W^{2,p}(I_L,\RN)$ and, in addition, that the left limit of $\dot\gc$ at $s=0$ agrees with the {right} limit at $s=L$, so that by the H\"older continuity of $\dot\gc$ we obtain condition
$ \dot\gc(0)=\dot\gc(L)$. In fact, if $\dot\gc(0)\neq\dot\gc(L)$ we get $\F_p(\gc)=+\i$ for any $p>1$.
\par
In a similar way, arguing as in Theorem~\ref{Tub} we obtain that if $\gc$ belongs to $W^{2,p}(I_L,\RN)$ for some $p>1$ and
$\dot\gc(0)=\dot\gc(L)$, then this time we can find a sequence $\{P_h\}$ of closed polygonal curves inscribed in $\gc$ such that $\m_\gc(P_h)\to 0$ and inequality \eqref{ub} holds. Therefore, the assertion readily follows.
\epf
\section{On the definition of p-rotation}\label{Subs:rmks}
In this final section, we briefly discuss whether alternative definitions of $p$-rotation of inscribed polygonals yield to a notion of $p$-energy for which
Theorems~\ref{Tlb} and \ref{Tub} continue to hold.

Since in all reasonable situations we have in mind we shall obtain the same $p$-energy, we thus conclude that our choice seems to be the more fitting one, at least from the point of view of numerical applications.
\par Coming back to the beginning of Section~\ref{Subs:prot}, we recall that $t_{i}^-$ and $t_i^+$ are the points in the segments $\gv_i$ and $\gv_{i+1}$, respectively,
whose distance to $v_i$ is equal to $r_i/2$, where $r_i:=\min\{\ell_i,\,\ell_{i+1}\}$ and $\ell_i=|\gv_i|$. Also, when the turning angle $\theta_i$ at $v_i$ is {in $]0,\pi[$,}  we chose by $\g_i$ the circular arc with end points $t_{i}^\pm$ and with tangent parallel to $\gv_i$ and
$\gv_{i+1}$ at $t_{i}^-$ and $t_i^+$, respectively (refer to Figure \ref{fig:pol} on page \pageref{fig:pol}).
On account of Definition~\ref{Dprot}, its {$p$-energy}, that is equal to $(r_i/2)^{1-p}\,\theta_i\,(\tan(\theta_i/2))^{p-1}$, gives the local
contribution to the $p$-rotation of the polygonal.
\par By a scaling argument, we e.g. may alternatively choose a term of the form $r_i^{1-p}\,f_p(\theta_i)$ for some angle function
$\theta\mapsto f_p(\theta)$ depending on the exponent $p$.
This yields to a different notion of $p$-rotation w.r.t. our choice for $\gk_p(P)$, where $f_p(\theta)=\theta\,(\tan(\theta/2))^{p-1}$.
For example, with $f^*_p(\theta):=2^{1-p}\theta^p$ one obtains a discrete $p$-curvature $\gk_p^*(P)$ that for equilateral polygonals agrees with formula
\eqref{Kstar}.
\par On the one hand, Theorem~\ref{Tub} continues to hold if
$$ \lim_{\theta\to 0^+}\frac{f_p(\theta)}{\theta\,(\tan(\theta/2))^{p-1}}=1 $$
that is the case of e.g. $f^*_p(\theta)=2^{p-1}\theta^p$. More generally, one can take
$$ f_{p,\a}(\t):=2^\a\,\t^{1-\a}\cdot\tan(\t/2)^{p-1+\a}\,,\qquad \a\geq 1-p $$
so that for $\a=1-p$ one gets $f^*_p(\theta)$ and for $\a=0$ our energy density for $\gk_p(P)$.
\par On the other hand, if one wishes that the assertion in Theorem~\ref{Tlb} continues to hold, one needs that $f_p(\theta)$ is equal to $\int_\gc \kk^p\,ds$, where $\kk$ is the scalar curvature of a curve $\gc$ with total curvature $\theta$ and satisfying the same first order boundary condition as the ones of a circular arc of radius one and total curvature $\theta$.
Actually, by Proposition~\ref{Pmin1}, with $\ell/2=1$, it turns out that this is the case for $f_{p,\a}(\t)$ provided that $\a\geq 0$, but not when $\a<0$, as e.g. for $f^*_p(\theta)$.
In fact, for $1-p\leq\a<0$, it is possible to find a curve $\gc$ satisfying the given clamping conditions and such that
$\int_\gc \kk^p\,ds=(\ell/2)^{p-1}\cdot f_{p,\a}(\t)$, by describing a great arc $\gc$ with very small curvature radius, but in that case the total curvature of $\gc$ would be
at least $2\pi-\t$, and hence greater than the turning angle of the polygonal $\theta$, when $\theta\in]0,\pi[$.
\par
As a consequence, the argument in the proof of Theorem~\ref{Tlb} fails to hold when $1-p\leq\a<0$.
In fact, on account of the lower semicontinuity inequality \eqref{lsc}, its validity depends on the existence of curves near $\gc$ whose distance goes to zero and whose $p$-energy is defined in correspondence to $f_{p,\a}(\t)$.
\par Therefore, denoting by $\gk_{p,\a}(P)$ the $p$-rotation of a polygonal $P$ in terms of the density $f_{p,\a}(\t)$, and letting
$$
 \F_{p,\a}(\gc):=\inf\Bigl\{\liminf_{h\to\i}\gk_{p,\a}(P_h)\mid \{P_h\}\ll\gc\,,\,\,\m_\gc(P_h)\to 0\Bigr\} $$
according to Definition~\ref{DEp}, it turns out that when $\a\geq 0$, both Theorems~\ref{Tlb} and \ref{Tub} continue to hold, and we actually obtain the same $p$-energy functional.
In fact, similarly to Theorem~\ref{Topen} we obtain:
\bc Let $\gc$ be a rectifiable and open curve in $\RN$ parameterized in arc-length, and let $\a\geq 0$. Then for every exponent $p>1$
$$ \F_{p,\a}(\gc)<\i \iff \gc\in W^{2,p}(I_L,\RN) $$
and in that case
$$ \F_{p,\a}(\gc)=\int_0^{L}\Vert\ddot \gc(s)\Vert^p\,ds\,. $$ \ec
\par Of course, all previous choices are not the optimal one in terms of energy minimizing inscribed curves, see Proposition~\ref{Pmin2} and Remark~\ref{Ropt2}.
The optimal choice among curves with total curvature equal to $\theta_i$, and with first order conditions at the middle points of the
consecutive segments $\gv_i$ and $\gv_{i+1}$ concurring at the vertex $v_i$, 
would be obtained by taking as $\g_i$ the energy minimizer of $\int_\gc\kk^p\,ds$ under the clamping conditions.
\par Denoting by $\gk_p^{{\textrm{opt}}}(P)$ the {\em optimal} $p$-rotation of a polygonal $P$ in terms of the latter choice, and
$$
 \F_p^{{\textrm{opt}}}(\gc):=\inf\Bigl\{\liminf_{h\to\i}\gk_p^{{\textrm{opt}}}(P_h)\mid \{P_h\}\ll\gc\,,\,\,\m_\gc(P_h)\to 0\Bigr\} $$
it turns out that the lower semicontinuity argument in the proof of Theorem~\ref{Tlb} continues to hold, and we thus readily obtain:
\bc Let $\gc$ be a rectifiable and open curve in $\RN$ parameterized in arc-length such that $\F_p^{{\textrm{opt}}}(\gc)<\i$ for some $p>1$.
Then $\gc\in W^{2,p}(I_L,\RN)$ and
$$ \int_0^{L} \Vert \ddot\gc(s)\Vert^p\,ds\leq \F_p^{{\textrm{opt}}}(\gc)<\i\,. $$
\ec
Therefore, since $\F_p^{{\textrm{opt}}}(\gc)\leq \F_p(\gc)$, by Theorem~\ref{Tub} we readily infer:
\bt\label{Topenopt} Let $\gc$ be a rectifiable and open curve in $\RN$ parameterized in arc-length. Then for every exponent $p>1$
$$ \F_p^{{\textrm{opt}}}(\gc)<\i \iff \gc\in W^{2,p}(I_L,\RN) $$
and in that case
$$ \F_p^{{\textrm{opt}}}(\gc)=\int_0^{L}\Vert\ddot \gc(s)\Vert^p\,ds\,. $$ \et
\par Again, the $p$-curvature functional we obtain is exactly the same.
However, apart from the case when $\Vert\gv_i\Vert=\Vert\gv_{i+1}\Vert$, for which we refer to Proposition~\ref{Pmin1},
it is not known {in general} how to compute the optimal curve, and hence the exact value of $\gk_p^{{\textrm{opt}}}(P)$,
see Remark~\ref{Ropt}.
\end{document}